\documentclass[12pt, a4paper, leqno]{article}
\usepackage[margin=2.5cm]{geometry}
\usepackage[utf8]{inputenc} 
\usepackage{lmodern}
\usepackage[T1]{fontenc}
\usepackage[english]{babel}

\usepackage[runin]{abstract}
\usepackage{titling}

\usepackage{amsmath}
\usepackage{amsthm}
\usepackage{amsfonts}
\usepackage{amssymb}
\usepackage{bbm}
\usepackage{mathtools}
\usepackage{clrscode}
\usepackage{enumitem}
\usepackage{mdwlist}

\abslabeldelim{.}

\newcommand{\keywordsname}{Key words}
\makeatletter
\newcommand{\keywords}[1]{%
\def\thekeywords{#1}%
\begin{@bstr@ctlist}
\hspace*{\abstitleskip}{\abstractnamefont\keywordsname\@bslabeldelim}\abstracttextfont\
#1%
\par\end{@bstr@ctlist}
}
\makeatother

\newcommand{\subjclassname}{Mathematics subject classification}
\makeatletter
\newcommand{\subjclass}[2][2020]{%
\begin{@bstr@ctlist}
\hspace*{\abstitleskip}{\abstractnamefont\subjclassname\ (#1)\@bslabeldelim}\abstracttextfont\
#2%
\par\end{@bstr@ctlist}
}
\makeatother

\makeatletter
\def\and{
	\end{tabular}%
	and%
	\begin{tabular}[t]{c}}%
\makeatother

\makeatletter
\def\thanks#1{
\protected@xdef\@thanks{\@thanks
\protect\footnotetext[\the\c@footnote]{#1}}%
}
\makeatother

\makeatletter
\let\addresses\@empty      
\newcommand{\address}[2][]{\g@addto@macro\addresses{\address{#1}{#2}}}
\newcommand{\curraddr}[2][]{\g@addto@macro\addresses{\curraddr{#1}{#2}}}
\newcommand{\email}[2][]{\g@addto@macro\addresses{\email{#1}{#2}}}
\newcommand{\urladdr}[2][]{\g@addto@macro\addresses{\urladdr{#1}{#2}}}
\def\enddoc@text{
  \ifx\@empty\addresses \else\@setaddresses\fi}
\AtEndDocument{\enddoc@text}
\def\emailaddrname{E-mail address}
\def\@setaddresses{\par
  \nobreak \begingroup
  \interlinepenalty\@M
  \def\address##1##2{\begingroup%
    \par\addvspace\bigskipamount
    \@ifnotempty{##1}{(\ignorespaces##1\unskip) }%
    {\noindent\ignorespaces##2}\par\endgroup}%
  \def\email##1##2{\begingroup
    \@ifnotempty{##2}{\nobreak\noindent\emailaddrname
      \@ifnotempty{##1}{, \ignorespaces##1\unskip}\/:\space
      \ttfamily##2\par}\endgroup}%
  \addresses
  \endgroup
}

\makeatother

\makeatletter
\def\cstar#1{\expandafter\@cstar\csname c@#1\endcsname}
\def\@cstar#1{\ifcase#1\or $\ast$\or $\ast\ast$\or $\ast\ast\ast$\fi}
\AddEnumerateCounter{\cstar}{\@cstar}{$\ast\ast\ast$}
\makeatother

\newcommand{\bfperiod}[1]{\textbf{#1.}}
\setlist[description]{leftmargin=0pt,format=\normalfont\bfperiod}

\newlist{conditions}{enumerate}{1}
\newlist{iconditions}{enumerate}{1}
\newlist{inthm}{enumerate}{1}

\setlist[conditions]{label=\normalfont(\alph*),ref=\normalfont(\alph*)}
\setlist[iconditions]{label=\normalfont(\roman*),ref=\normalfont(\roman*)}
\setlist[inthm]{label=\normalfont(\thetheorem.\arabic*),ref=\normalfont(\thetheorem.\arabic*),wide,labelindent=0pt}

\mathchardef\mhyphen="2D

\newcommand{\R}{\mathbb{R}}

\newcommand{\C}{\mathcal{C}}

\newtheorem{theorem}{Theorem}[section]
\newtheorem{corollary}[theorem]{Corollary}
\newtheorem{lemma}[theorem]{Lemma}
\newtheorem{fact}[theorem]{Fact}
\newtheorem{proposition}[theorem]{Proposition}
\theoremstyle{definition}

\newtheorem{example}[theorem]{Example}
\newtheorem{remark}[theorem]{Remark}

\newtheorem{df}[theorem]{Definition}

\newtheorem{claim}[theorem]{Claim}
\newtheorem{observation}[theorem]{Observation}

\DeclarePairedDelimiter\abs{\lvert}{\rvert}%
\DeclarePairedDelimiter\norm{\lVert}{\rVert}%

\makeatletter
\let\oldabs\abs
\def\abs{\@ifstar{\oldabs}{\oldabs*}}
\let\oldnorm\norm
\def\norm{\@ifstar{\oldnorm}{\oldnorm*}}
\makeatother

\makeatletter
\renewcommand\paragraph{
  \@startsection{paragraph}
                {4}
                {\z@}
                {3.25ex \@plus1ex \@minus.2ex}
                {-0.4em}
                {\normalfont\normalsize\bfperiod}}
\makeatother

\title{Approximation of maps from algebraic polyhedra to real algebraic varieties}
\date{}
\author{Marcin Bilski \and Wojciech Kucharz}

\address{Marcin Bilski\\Institute of Mathematics\\Faculty of Mathematics and Computer
Science\\Jagiellonian University\\\L{}ojasiewicza 6\\30-348
Krak\'ow\\Poland}
\email{Marcin.Bilski@im.uj.edu.pl}

\address{Wojciech Kucharz\\Institute of Mathematics\\Faculty of Mathematics and Computer
Science\\Jagiellonian University\\\L{}ojasiewicza 6\\30-348
Krak\'ow\\Poland}
\email{Wojciech.Kucharz@im.uj.edu.pl}

\usepackage[pdftex, pdfauthor={M. Bilski, J. Bochnak, W. Kucharz}, pdftitle={\thetitle}]{hyperref}

\begin{document}
\maketitle
\thispagestyle{empty}

\begin{abstract}
Given a finite simplicial complex $\mathcal{K}$ in $\R^n$ and a real algebraic variety $Y,$ by a \linebreak$\mathcal{K}$-regular map $|\mathcal{K}|\rightarrow Y$ we mean a continuous map whose restriction to every simplex in $\mathcal{K}$ is a regular map. A simplified version of our main result says that if $Y$ is a uniformly retract  rational variety and if $k, l$
are integers satisfying $0\leq l\leq k,$ then every $\mathcal{C}^l$ map $|\mathcal{K}|\rightarrow Y$ can be approximated in the $\mathcal{C}^l$ topology by $\mathcal{K}$-regular maps of class $\mathcal{C}^k.$ By definition, $Y$ is uniformly retract rational if for every point $y\in Y$ there is a Zariski open neighborhood $V\subset Y$ of $y$ such that the identity map of $V$ is the composite of regular maps $V\rightarrow W\rightarrow V,$ where $W\subset\R^p$ is a Zariski open set for some $p$ depending on $y.$  
\end{abstract}


\phantomsection
\addcontentsline{toc}{section}{\refname}


\section{Introduction}
\label{intro}

Throughout this paper we use the term \textit{real algebraic variety} to mean a ringed space with structure sheaf
of $\mathbb{R}$-algebras  of $\mathbb{R}$-valued functions, which is isomorphic to a Zariski locally closed set
(in some real projective space) endowed with the Zariski topology and the sheaf of regular functions. This is
compatible with the monographs  \cite{bib4, ManBook} containing a detailed exposition of real algebraic geometry.
Morphisms of real algebraic varieties are called \textit{regular maps.}

Recall that each real algebraic variety in the sense used here is actually \textit{affine,} that is, isomorphic to
an algebraic set in $\mathbb{R}^n$ for some $n$ \cite[Proposition 3.2.10 and Theorem 3.4.4]{bib4} or
\cite[Proposition 1.3.11]{ManBook}. This fact allows us to describe regular functions and regular maps in a simple
way. To be more precise, let $X\subset\mathbb{R}^n$ and $Y\subset\mathbb{R}^m$ be algebraic sets, and let
$U\subset X$ be a Zariski open set. A function $\varphi:U\rightarrow\mathbb{R}$ is regular if and only if there exist
two polynomial functions $P, Q:\mathbb{R}^n\rightarrow\mathbb{R}$ such that
$$ Q(x)\neq 0\mbox{ and }\varphi(x)=\frac{P(x)}{Q(x)} \mbox{ for all  }x\in U$$
\cite[p. 62]{bib4} or \cite[p. 14]{ManBook}. A map $f=(f_1,\ldots,f_m):U\rightarrow Y\subset\mathbb{R}^m$ is
regular if and only if its components $f_1,\ldots,f_m$ are regular functions.

Each real algebraic variety is also equipped with the Euclidean topology determined by the usual metric on
$\mathbb{R}.$ Unless explicitly stated otherwise, all topological notions relating to real algebraic varieties
will refer to the Euclidean topology.


According to the classical Weierstrass theorem, if $S\subset\mathbb{R}^n$ is  a compact set, then each continuous
real-valued function on $S$ can be approximated uniformly by polynomial  functions. However, as indicated in
Example \ref{ex-2-2} below, polynomial maps are too rigid for approximation of maps with values in real algebraic varieties (even so simple varieties as the Euclidean unit sphere in $\R^m$).

The problem of approximation of continuous maps between real algebraic varieties by regular maps or by more
flexible "almost regular" maps has been investigated in numerous works (see \cite{Ban, BanKu, BenWit, bib2, bib3, bib4,
bib6, bib7, bib77, Ku2009, Ku2014, bib18, bib19, bib20, bib21, ManBook, Zie} and the references therein). "Almost
regular" is interpreted in the literature to mean one of the following: (1) continuous rational \cite{Ku2009,
Ku2014}, (2) regulous \cite{bib13, bib16, bib17, bib19, bib20, bib21, Zie} (regulous=hereditarily
rational=stratified-regular, and for nonsingular varieties, regulous=continuous rational), (3) quasi-regulous
\cite{bib2} and (4) piecewise-regular \cite{bib3, bib18}. Approximation properties of regular maps and almost
regular maps are briefly reviewed in \cite{bib3}. Fernando and Ghiloni \cite{FeGh} give an overview of other
aspects of the approximation problem and prove remarkable theorems on approximation by differentiable
semialgebraic maps, which are further generalized by Paw\l ucki \cite{Paw} and Carbone \cite{Carb}. In \cite{BenWit}, Benoist and Wittenberg obtain a deep result on regular approximation for
maps defined on nonsingular real algebraic curves. The recent papers \cite{Ban,BanKu,bib77} and \cite{bib19} present the
state of the art of approximation by regular maps and by regulous maps, respectively.

In the present paper, we consider approximation by maps constituting a proper subclass of the class of piecewise-regular ones.  Before we define approximating maps, let us recall the generalization of the notion of regular map introduced in
\cite{bib18}.
\begin{df}\label{def1}Let $X, Y$ be real algebraic varieties, $S$ an arbitrary subset of $X,$ $T$ an arbitrary subset of $Y,$ and $Z$ the Zariski
closure of $S$ in $X.$ A map $f:S\rightarrow T$ is said to be \textit{regular} if there exist a Zariski open neighborhood
$Z_0\subset Z$ of $S$ and a regular map (in the usual sense) $f_0:Z_0\rightarrow Y$ such that $f_0(x)=f(x)$ for all $x\in S.$
\end{df}

In this general context, the composite of regular maps is a regular map, so the notion of biregular isomorphism makes sense.   

By a \textit{simplex} in $\mathbb{R}^n$ we always mean a closed geometric simplex. A subset $C$ of $\mathbb{R}^n$
is called an \textit{algebraic cell} of dimension $d$ if there exists a biregular isomorphism $\varphi:\Delta\rightarrow C$ from
a $d$-dimensional simplex $\Delta$ in $\R^n$ onto $C.$ An $e$-dimensional face of $C$ is, by definition, the image of an
$e$-dimensional face of $\Delta$ by $\varphi.$ Clearly, the notion of face is independent of the choice of biregular
isomorphism $\varphi:\Delta\rightarrow C.$ It is also clear that an $e$-dimensional face of $C$ is an algebraic cell of
dimension $e.$

An \textit{algebraic complex} $\mathcal{K}$ in $\mathbb{R}^n$ is a finite collection of
algebraic cells in $\R^n$ such that: (i) if $K\in\mathcal{K},$ then every face of $K$ is in
$\mathcal{K}$ and (ii) if $K\cap L\neq\emptyset$ for some $K, L\in\mathcal{K},$ then $K\cap L$ is the union of some faces of
both $K$ and $L.$ For any algebraic complex $\mathcal{K},$ we write $|\mathcal{K}|$ to denote the \textit{algebraic
polyhedron} defined to be the union of all algebraic cells in $\mathcal{K}.$ Clearly, $|\mathcal{K}|$ is a compact
semialgebraic subset of $\mathbb{R}^n.$
\begin{df}\label{def2}Let $\mathcal{K}$ be an algebraic complex in $\mathbb{R}^n$ and let $Y$ be a real algebraic
variety. A map $f:|\mathcal{K}|\rightarrow Y$ is said to be \textit{$\mathcal{K}$-regular} if for every algebraic cell
$C\in\mathcal{K}$ the restriction $f|_C:C\rightarrow Y$ is a regular map (in particular $f$ is continuous). More
generally, a map $g:S\rightarrow Y$ defined on an arbitrary subset $S$ of $|\mathcal{K}|$ is said to be
$\mathcal{K}$-regular if it is the restriction of some $\mathcal{K}$-regular map from $|\mathcal{K}|$ to $Y.$
\end{df}


We will consider approximation by $\mathcal{K}$-regular maps of class $\mathcal{C}^k$ 
in the space of $\mathcal{C}^l$ maps for $0\leq l\leq k<+\infty.$

First recall what "class $\mathcal{C}^l$" means in case of maps defined on subsets of $\R^n$ which do not have differentiable structure. 
Let $S\subset\R^n$ be a compact set, $Y$ a nonsingular real algebraic variety, and 
$l$ a nonnegative integer. A map $f:S\rightarrow Y$ is said to be of \textit{class} $\mathcal{C}^l,$ or a $\mathcal{C}^l$ \textit{map}, if it extends to a $\mathcal{C}^l$ {map} (in the usual sense) $U\rightarrow Y$ defined on an open neighborhood $U\subset\R^n$ of $S.$ Denote by ${\mathcal{C}}^l(S,Y)$ the set of all $\mathcal{C}^l$ maps from $S$ to $Y.$

In \cite[Section 2.3]{bib3} we defined a natural topology on ${\mathcal{C}}^l(S,Y)$, called the ${\mathcal{C}}^l$ \textit{topology}. Here we only point
out that the $\mathcal{C}^0$ topology on ${\mathcal{C}}^0(S,Y)$ is the compact-open topology, and the $\mathcal{C}^l$ topology on $\mathcal{C}^l(S,Y)$ is the compact-open  $\mathcal{C}^l$ topology  (discussed in \cite[pp. 34, 35]{Hi})
if $S$ is a compact $\mathcal{C}^{\infty}$ submanifold of $\R^n.$

Let $f\in\mathcal{C}^l(S,Y)$ and let $\mathcal{A}$ be a subset of
$\mathcal{C}^k(S,Y),$ where $k,l$ are integers satisfying $0\leq l\leq k.$ 
We say that $f$ can be \textit{approximated in the $\mathcal{C}^l$ topology  by maps from} $\mathcal{A}$ if every open neighborhood $\mathcal{V}$ of $f$ in $\mathcal{C}^l(S,Y)$ contains a map
from $\mathcal{A}$. This is equivalent to the following: Fix an algebraic embedding $Y\subset \R^m.$ 
If $\varphi=(\varphi_1,\ldots,\varphi_m):U\rightarrow\R^m$ is an extension
of class $\mathcal{C}^l$ of $f$ defined on an open neighborhood $U\subset\R^n$ of $S,$ 
then for every constant $\varepsilon>0$ there exists a $\mathcal{C}^k$ map 
$\psi=(\psi_1,\ldots,\psi_m):U\rightarrow\R^m$ such that $\psi(S)\subset Y,$ 
the restriction $\psi|_S:S\rightarrow Y$ belongs to $\mathcal{A},$ and

$$\big{|}\frac{\partial^{|\alpha|} {\psi_i}}{\partial x_1^{\alpha_1}\cdots\partial x_n^{\alpha_n}}(x)-\frac{\partial^{|\alpha|} {\varphi_i}}{\partial x_1^{\alpha_1}\cdots\partial x_n^{\alpha_n}}(x)\big{|}<\varepsilon$$ for all $x\in S,$ $1\leq i\leq m,$ and $\alpha=(\alpha_1,\ldots,\alpha_n)\in\mathbb{N}^n$ with $|\alpha|:=\alpha_1+\ldots+\alpha_n\leq l$ ($\mathbb{N}=\{0,1,\ldots\}$).

Given an algebraic complex $\mathcal{K}$ in $\R^n,$ it is easy to see that the $\mathcal{K}$-regular maps
$|\mathcal{K}|\rightarrow Y$ constitute a proper subclass of the class of piecewise-regular maps from
$|\mathcal{K}|$ to $Y$ defined in \cite{bib18}. By \cite[Theorem 1.1]{bib3}, for any nonnegative integers $l\leq
k,$ every $\mathcal{C}^l$ map from $|\mathcal{K}|$ to $Y$ can be approximated 
in the $\mathcal{C}^l$ topology by  piecewise-regular maps of class $\mathcal{C}^k$ under the assumption that $Y$ is a \textit{uniformly rational variety} (the latter notion is recalled in Section~\ref{Sec2} below).  
In the present paper, we will prove a stronger result allowing approximation 
by $\mathcal{K}$-regular maps 
of class $\mathcal{C}^k$ with a weaker assumption on the target variety $Y$; the corresponding class of target varieties is defined immediately below.

\begin{df}\label{def3} A real algebraic variety $Y$ is said to be \textit{uniformly retract rational} if for every point $y\in Y$ there exist a Zariski open neighborhood $V$ of $y$ in $Y,$ a Zariski open subset $W$ of $\R^p$ (for some $p$ depending on $y$), and two regular maps $i:V\rightarrow W,$ $r:W\rightarrow V$ such that the composite $r\circ i:V\rightarrow V$ is the identity map.
\end{df}
This definition, introduced by Banecki \cite{Ban}, is a variant of the standard concept of
retract rational variety that has been present in the literature since the 1980s. By \cite[Observation 1.4]{Ban}, every uniformly retract rational real algebraic variety is nonsingular.


What follows is our main result.
\begin{theorem}\label{mainr} Let $\mathcal{K}$ be an algebraic complex in $\mathbb{R}^n$ and let $Y$ be a uniformly retract rational
real algebraic variety. Let $k, l$ be
integers satisfying $0\leq l\leq k$. Then every $\mathcal{C}^l$ map from $|\mathcal{K}|$ to $Y$ can be approximated in the
$\mathcal{C}^l$ topology by $\mathcal{K}$-regular maps $|\mathcal{K}|\rightarrow Y$ of class $\mathcal{C}^k.$
\end{theorem}

For source spaces that are not necessarily algebraic complexes, Theorem~\ref{mainr} leads to the following result.

\begin{corollary}\label{corofmainr} Let $S\subset\mathbb{R}^n$ be a locally contractible compact set
and let $Y$ be a uniformly retract rational real algebraic variety. Let $k, l$ be integers satisfying $0\leq l\leq k.$ Then there exists a finite simplicial complex $\mathcal{K}$ in $\mathbb{R}^n$ with $S\subset |\mathcal{K}|$ such that every $\mathcal{C}^l$ map from $S$ to $Y$ can be
approximated in the $\mathcal{C}^l$ topology by $\mathcal{K}$-regular maps $S\rightarrow Y$ of class $\mathcal{C}^k.$ 
\end{corollary}
\proof   By Borsuk's theorem \cite[p. 537, Theorem E.3]{bib9}, there exists a continuous retraction $r:U\rightarrow S$ of a neighborhood $U$ of $S$ in $\mathbb{R}^n.$ Since $S$ is a compact subset of $U,$ we obtain $S\subset B\subset U,$ where $B$ is the union of a finite collection of simplices in $\mathbb{R}^n.$ Clearly, there exists a finite simplicial complex $\mathcal{K}$ in $\mathbb{R}^n$ with $|\mathcal{K}|=B.$   

Let $f:S\rightarrow Y$ be a $\mathcal{C}^l$ map. Then there is
a $\mathcal{C}^l$ extension $\tilde{f}:V\rightarrow Y$ of $f$ to some open neighborhood $V$ of $S$ in $\mathbb{R}^n.$   By approximating $r$ relatively to $S$ we obtain a $\mathcal{C}^{\infty}$ map $g:U\rightarrow V$ such that $g|_S$ is the identity on $S.$ Then the composite $\tilde{f}\circ g|_{|\mathcal{K}|}:|\mathcal{K}|\rightarrow Y$ is a $\mathcal{C}^l$ extension of $f.$ By Theorem~\ref{mainr}, $\tilde{f}\circ g|_{|\mathcal{K}|}$ can be approximated in the $\mathcal{C}^l$ topology by $\mathcal{K}$-regular maps of class $\mathcal{C}^k$, which completes the proof. 
\qed\vspace*{2mm}

In a geometric setting, local contractibility is only a minimally restrictive condition. For example, every
semialgebraic subset of $\mathbb{R}^n$ is locally contractible \cite[Corollary 9.3.7]{bib4}.


Let us note that Theorem~\ref{mainr} and Corollary~\ref{corofmainr} do not hold in general for $k=\infty$ (cf. \cite[Example 1.7]{bib18}). Moreover, the assumptions on the target variety $Y$ cannot be relaxed too much. To see this, recall that a real algebraic variety $Y$ is said to be \textit{unirational} if there exists a regular map $\varphi:U\rightarrow Y$ defined on a Zariski open subset $U$ of $\R^n,$ for some $n,$ such that the image $\varphi(U)$ is Zariski dense in $Y.$ If $Y$ is unirational, then it is irreducible and we can choose a regular map $\varphi$ as above with $n=\mathrm{dim}Y.$

As usual, for any nonnegative integer $n$, let $\mathbb{S}^n$ denote the unit $n$-sphere, $$\mathbb{S}^n:=\{(x_0,\ldots,x_n)\in\mathbb{R}^{n+1}:x_0^2+x_1^2+\cdots+x_n^2=1\}.$$  

\begin{observation}\label{obs-1-4}
Let $Y$ be a $p$-dimensional irreducible real algebraic variety and let $\Delta$ be a $d$-simplex in $\R^n$ with $d\geq p\geq 1.$ Assume that every continuous map from $\Delta$ to $Y$ can be approximated 
in the $\mathcal{C}^0$ topology by regular maps. We claim that then $Y$ is a unirational variety.

In the rest of this observation, approximation will mean approximation in the compact-open topology. If $\Gamma$ is a $p$-dimensional face of $\Delta,$ then there is a retraction from $\Delta$ onto $\Gamma,$ and hence every continuous
map from $\Gamma$ to $Y$ can be approximated by regular maps. Thus, the claim is reduced to the case $d=p.$ Moreover, we may assume without loss of generality that $Y$ is a nonsingular variety. It suffices to construct a continuous map $f:\Delta\rightarrow Y$ such that for each continuous map $g:\Delta\rightarrow Y$ sufficiently close to $f$ the image $g(\Delta)$ contains a nonempty open subset of $Y.$
Indeed, in this case, if $g$ is a regular map close to $f,$ then $g(\Delta)$ is Zariski dense in $Y,$ so $Y$ is a unirational variety.
The construction of $f$ is a purely topological task.

For any constant $r>0,$ let $\overline{B_r}$ denote the closed Euclidean ball in $\R^p$ centered at the origin with radius $r.$ We will work with $\overline{B_1}$ instead of $\Delta,$ which is allowed because these two spaces are homeomorphic. 
Choose an open subset $V$ of $Y$ and a homeomorphism $\alpha:\R^p\rightarrow V.$ It remains to show that the restriction
$f:\overline{B_1}\rightarrow V\subset Y$ of $\alpha$ has the asserted approximation property. 
If a continuous map $g:\overline{B_1}\rightarrow Y$ is sufficiently close to $f,$ then $g(\overline{B_1})\subset V,$ so composing with $\alpha^{-1}$ reduces the task to the case where $f:\overline{B_1}\rightarrow\R^p$ is the inclusion map and $g:\overline{B_1}\rightarrow\R^p$ is close to $f.$ Thus, to complete the argument, it suffices to prove the following:    
if $h:\overline{B_1}\rightarrow\R^p$ is a continuous map satisfying $$||h(x)-x||<\frac{1}{2}\mbox{ for all }x\in \overline{B_1},$$ then $\overline{B_{\frac{1}{2}}}\subset h(\overline{B_1})$ (here $||\cdot||$ stands for the Euclidean norm). Suppose to the contrary that there is a point $z\in \overline{B_{\frac{1}{2}}}\setminus h(\overline{B_1}).$ Define a continuous map $H:\overline{B_1}\rightarrow\mathbb{S}^{p-1}$ by sending $x\in \overline{B_1}$ to the intersection of $\mathbb{S}^{p-1}$ with the ray emanating from $z$ through $h(x).$ It follows that the restriction $H|_{\mathbb{S}^{p-1}}:\mathbb{S}^{p-1}\rightarrow\mathbb{S}^{p-1}$ is a null homotopic map. On the other hand, by construction, $H(x)\neq -x$ for all $x\in\mathbb{S}^{p-1},$ so the map $H|_{\mathbb{S}^{p-1}}$ is homotopic to the identity map of $\mathbb{S}^{p-1}$ via the homotopy $$\mathbb{S}^{p-1}\times[0,1]\rightarrow\mathbb{S}^{p-1}, \mbox{ } (x,t)\mapsto\frac{(1-t)H(x)+tx}{||(1-t)H(x)+tx||},$$ a contradiction.   
\end{observation}

There is  a good reason why in this paper we work with maps whose restrictions to simplices are regular rather than polynomial maps. By definition, a \textit{polynomial map} from a subset of $\mathbb{R}^n$ to a subset of $\mathbb{R}^m$ is the restriction of a polynomial map $\R^n\rightarrow\R^m.$

\begin{example}\label{ex-2-2}Let $d, n, p, q$ be positive integers, $\Delta$ a $d$-simplex in $\R^n,$ and $f:\Delta\rightarrow\mathbb{S}^p$  a continuous map. Assume that $d>q>p$ and $q$ is an integer power of $2.$ Then the following conditions are equivalent:   \vspace*{3mm}\\
(a) $f$ can be approximated in the $\mathcal{C}^0$ topology by polynomial maps $\Delta\rightarrow\mathbb{S}^p.$\\
(b) $f$ is a constant map.\vspace*{3mm}\\
Suppose that $f$ is not a constant map. We can choose an affine-linear $(q+1)$-plane $H\subset\R^n$ and a $q$-dimensional sphere $S\subset H$ such that $S\subset H\cap\Delta$ and the restriction $f|_S$ is not constant. Let $\varphi:\mathbb{S}^q\rightarrow S$ be the restriction of an affine-linear isomorphism $\R^{q+1}\rightarrow H.$ Then the composite map $f\circ\varphi:\mathbb{S}^q\rightarrow\mathbb{S}^p$ is not constant.  Thus, if there were a polynomial map $g:\Delta\rightarrow\mathbb{S}^p$ close enough to $f,$ the composite polynomial map $g\circ\varphi:\mathbb{S}^q\rightarrow\mathbb{S}^p$ would not be constant, contradicting Wood's theorem
\cite{Wo} or \cite[Theorem 13.1.9]{bib4}. Consequently, conditions (a) and (b) are equivalent.
\end{example}

The paper is organized as follows. In Section~\ref{Sec2} we recall the notion of uniformly rational variety and  discuss examples showing interesting instances of applicability of Theorem~\ref{mainr}. 
In Section \ref{sec-3}  we develop some tools which, together with Banecki's recent paper \cite{Ban}, allow us to prove our main result in Section \ref{sec-4}.

\section{Uniformly rational varieties}
\label{Sec2}

An $n$-dimensional irreducible real algebraic variety is said to be \textit{rational} if it contains a nonempty Zariski
open subset biregularly isomorphic to a Zariski open subset of $\mathbb{R}^n.$ 
An $n$-dimensional real algebraic variety (possibly reducible) is said to be \textit{uniformly rational} if each of its points has a Zariski
open neighborhood that is biregularly isomorphic to a Zariski open subset of $\mathbb{R}^n.$

Clearly, every uniformly rational real algebraic variety is nonsingular of pure dimension. It is an open question whether  every irreducible nonsingular rational
variety is uniformly rational (see \cite{BoBo} and \cite[p. 885]{Gro}
for the discussion involving complex algebraic varieties). 
\begin{observation} The relevance of the notion of uniformly rational variety to Theorem~\ref{mainr} stems from the following obvious fact: if a real algebraic variety is uniformly rational, then it is also uniformly retract rational. It is known that the converse is not true in general. Indeed, the complex algebraic hypersurface $V\subset\mathbb{C}^4$ defined in \cite[p. 315, Example 3]{BeCoSaSw} is not rational. Since $V$ is defined over $\R,$ the real algebraic hypersurface $X:=V\cap \R^4\subset\R^4$ is not rational. On the other hand, by \cite[p. 299, Th\'eor\`eme 1']{BeCoSaSw}, there exists a nonempty Zariski open subset of $X\times\R^3$ that is biregularly isomorphic to a Zariski open subset of $\R^6.$ Thus, some nonempty Zariski open
subset $Y$ of $X$ is a uniformly retract rational real algebraic variety. Obviously, $Y$ is not uniformly rational. We are indebted to Olivier Benoist for bringing the relevant results of \cite{BeCoSaSw} to our attention. 
\end{observation}

A \textit{linear real algebraic group} $G$ is a Zariski closed subgroup of the general linear group $\mathrm{GL}_n(\mathbb{R}),$ for
some $n.$ A \textit{$G$-space} is a real algebraic variety $Y$ on which $G$ acts, the action $G\times Y\rightarrow Y,$
$(a,y)\mapsto a\cdot y$ being a regular map. If the action is transitive,  then $Y$ is called a \textit{homogeneous space}
for $G.$ The group $G$ acts transitively on itself by left multiplication and is therefore a homogeneous space. If a homogeneous space has at least one rational irreducible component, then it is a uniformly rational
variety. It is an open question whether every linear real algebraic group is a unirational variety, but there are
homogeneous spaces for some linear real algebraic groups that are not uniformly rational varieties.
 
In the rest of this section, we will focus on uniformly rational varieties.    
The following example sheds some light on the scope of applicability of Theorem~\ref{mainr}.

\begin{example}\label{ex-2-1}Here are some uniformly rational real algebraic varieties.

\begin{inthm}[widest=1.3.3]
\item\label{ex-2-1-1}\textbf{(Unit spheres).} For any nonnegative integer $n$ the unit $n$-sphere $\mathbb{S}^n$ is a uniformly rational variety because $\mathbb{S}^n$ with any point removed is biregularly isomorphic to $\mathbb{R}^n$ via a suitable stereographic projection.

\item\label{ex-2-1-2}\textbf{(Grassmannians).} Let $\mathbb{F}$ stand for $\R, \mathbb{C}$ or $\mathbb{H},$ where $\mathbb{H}$ is the skew field of quaternions. Denote by $\mathbb{G}_r(\mathbb{F}^n)$ the Grassmannian of $r$-dimensional $\mathbb{F}$-vector subspaces of $\mathbb{F}^n,$ with $0\leq r\leq n,$ considered as  a real algebraic variety \cite[pp. 72, 73, 352]{bib4}. Clearly, $\mathbb{G}_r(\mathbb{F}^n)$ is a uniformly rational variety.
 
 \item\label{ex-2-1-3}\textbf{(Special orthogonal groups).} The special orthogonal group $\mathrm{SO}_n(\mathbb{R})\subset \mathrm{GL}_n(\mathbb{R})$ is a rational variety, so it is uniformly rational. Indeed, let $V_n$ be the space of all skew-symmetric $n$-by-$n$ matrices with real entries; thus $V_n$ is isomorphic to $\R^d,$ where $d={n(n-1)}/{2}.$ If $I\in \mathrm{SO}_n(\mathbb{R})$ is the identity matrix and $A\in V_n,$ then $I+A$ is an invertible matrix, and the Cayley transform
 $$V_n\rightarrow\{Q\in \mathrm{SO}_n(\mathbb{R}):\mathrm{det}(I+Q)\neq 0\}, A\mapsto(I-A)(I+A)^{-1} $$ is a biregular isomorphism.

 \item\label{ex-2-1-4}\textbf{(Orthogonal groups).} The orthogonal group $\mathrm{O}_n(\R)\subset \mathrm{GL}_n(\R)$ is a uniformly rational variety because $\mathrm{SO}_n(\R)$ is an irreducible component of $\mathrm{O}_n(\R).$
 
 \item\label{ex-2-1-5}\textbf{(Real form of complex groups).} By Chevalley's theorem \cite[Corollary 2]{bib10}, every linear complex algebraic group $G\subset\mathrm{GL}_n(\mathbb{C})$ has the property that each of its irreducible components is a rational variety. Thus, the image of $G$ under the standard embedding $\mathrm{GL}_n(\mathbb{C})\rightarrow \mathrm{GL}_{2n}(\R)$ is a linear real algebraic group that is a uniformly rational variety.   

 \item\label{ex-2-1-6}\textbf{(Blow-ups of uniformly rational varieties).} If $Y$ is a uniformly rational real algebraic variety and $Z$ is a nonsingular Zariski closed subvariety of $Y,$ then the blow-up of $Y$ with center $Z$ is a uniformly rational variety. The proof given in \cite{BoBo} and \cite[p. 885]{Gro} in a complex setting also works for real algebraic varieties.
 \end{inthm}

\end{example}

\section{Approximation with interpolation of functions}\label{sec-3}

Let $U\subset \R^n$ be an open set, $l$ a nonnegative integer, and $g:U\rightarrow\R$ a $\C^l$ function.   
We will often write $\partial^{\alpha}g$ to denote 
$\frac{\partial^{|\alpha|}g}{\partial x_1^{\alpha_1}\ldots\partial x_n^{\alpha_n}},$
for $\alpha\in\mathbb{N}^n$ with $|\alpha|\leq l.$
More generally, if $Y$ is a nonsingular Zariski closed subvariety of $\R^m$ and $h=(h_1,\ldots,h_m):U\rightarrow Y\subset\R^m$ is a $\C^l$ map, we set $$\partial^{\alpha}h:=(\partial^{\alpha}h_1,\ldots,\partial^{\alpha}h_m),$$
for $\alpha$ as above.

Let $S\subset\R^n$ be a compact set in $\R^n$ and let  
$$\varphi=(\varphi_1,\ldots,\varphi_m):V\rightarrow Y\subset\R^m,\mbox{ } 
\psi=(\psi_1,\ldots,\psi_m):W\rightarrow Y\subset\R^m$$ be $\C^l$ maps defined on  open 
neighborhoods $V, W$ of $S$ in $\R^n.$ Given a constant $\varepsilon>0,$ we say that $\varphi$ and $\psi$ are \textit{$\varepsilon$-close on $S$ in the $\mathcal{C}^l$ topology} if
$$|\partial^{\alpha}\varphi_i(x)-\partial^{\alpha}\psi_i(x)|<\varepsilon$$ for all $x\in S,$ $1\leq i\leq m,$ and $\alpha\in\mathbb{N}^n$ with $|\alpha|\leq l.$ Leaving $\varepsilon$ to be specified later, we can also simply say that $\varphi$ and $\psi$ are close on $S$ in the $\mathcal{C}^l$ topology.  

In this section we consider the approximation of real-valued functions, deferring the general case to the next section.

For any subset $A\subset\R^n$ let $\overline{A}^Z$  denote the closure of $A$ in the Zariski topology.

\begin{lemma}\label{cubeextension}Let $U\subset\R^n$ be an open
neighborhood of a $d$-simplex $\Delta$ in $\R^n,$ with $d\geq 1.$  Let $k$ be a positive integer and let $f:U\rightarrow\mathbb{R}$ be a $\mathcal{C}^{(d+1)k}$ function such that $\partial^{\alpha}f|_{\overline{\Gamma}^Z\cap U}$ is a regular function for every $(d-1)$-dimensional face $\Gamma$ of $\Delta$ and every $\alpha\in\mathbb{N}^n$ with $|\alpha|\leq k-1.$ Then for every constant $\varepsilon>0$ there exists a regular function $\tilde{f}:\tilde{U}\rightarrow\R$ defined on an open neighborhood $\tilde{U}\subset U$ of $\Delta$ such that $\tilde{f}$ and $f$ are $\varepsilon$-close on $\Delta$ in the $\mathcal{C}^k$ topology, and 
$$\partial^{\alpha}\tilde{f}|_{\overline{\Gamma}^Z\cap\tilde{U}}=\partial^{\alpha}f|_{\overline{\Gamma}^Z\cap\tilde{U}}$$ for every $(d-1)$-dimensional
face $\Gamma$ of $\Delta$ and every $\alpha\in\mathbb{N}^n$ with $|\alpha|\leq k-1.$
\end{lemma}

To prove Lemma~\ref{cubeextension}, we need the following fact, which seems to be well known, but we recall its proof for completeness.

\begin{fact}\label{ckdivision} Let $U$ be an open subset in $\mathbb{R}^n$ and let $f:U\rightarrow\mathbb{R}$ be a
$\mathcal{C}^k$ function, where $k$ is a positive integer. Let $H:\R^n\rightarrow\R$ be a nonzero affine-linear function. Assume that $f$ vanishes identically on
$H^{-1}(0)\cap U.$ Then the quotient
${f}/{H}:U\setminus H^{-1}(0)\rightarrow\mathbb{R}$ can be extended to a $\mathcal{C}^{k-1}$ function
$U\rightarrow\mathbb{R}.$
\end{fact}
\proof Without loss of generality we assume, applying an affine coordinate change if necessary, that $H(x_1,\ldots,x_n)=x_1.$ Choose an open subset $\tilde{U}$ of $U$
containing $H^{-1}(0)\cap U$ and such that for every point $(x_1,x_2,\ldots,x_n)$
in $\tilde{U}$ the point $(tx_1,x_2,\ldots,x_n)$ is also in $\tilde{U}$ for all $t\in[0,1].$ Then we have 
$$f(x_1,x_2,\ldots, x_n)=\int_0^1\frac{df(tx_1,x_2,\ldots,x_n)}{dt}dt=x_1\int_0^1\frac{\partial f
}{\partial x_1}(tx_1,x_2,\ldots,x_n)dt.$$ The proof is complete since the integral $\int_0^1\frac{\partial f }{\partial
x_1}(tx_1,x_2,\ldots,x_n)dt$ defines a function $\tilde{U}\rightarrow\mathbb{R}$ of class $\mathcal{C}^{k-1}$ which is an extension of the
restriction of ${f}/H$ to $\tilde{U}\setminus H^{-1}(0).$\qed\vspace*{4mm}\\
\textit{Proof of Lemma~\ref{cubeextension}.} Suppose first that $n=d.$ Let $H_{\Gamma}:\R^n\rightarrow\R$ be a nonzero affine-linear function vanishing on the face $\Gamma$ of $\Delta.$
Put $u(x_1,\ldots,x_n):=\prod_{\Gamma}H_{\Gamma}(x_1,\ldots,x_n),$ where $\Gamma$ in the product ranges over the set of all $(d-1)$-dimensional faces of $\Delta.$

\begin{claim}\label{claim03} There are an open
neighborhood $\tilde{U}$ of $\Delta$ in $U,$ a regular function $h:\tilde{U}\rightarrow\mathbb{R}$ and a
$\mathcal{C}^k$ function $g:\tilde{U}\rightarrow\mathbb{R}$ such that
$$f-h=gu^k \mbox{ on }\tilde{U}.$$
\end{claim}
\noindent\textit{Proof of Claim~\ref{claim03}.} For any integer $s$ with $1\leq s\leq k$ define:\vspace*{2mm}\\
\textit{Condition(s).} There are an open neighborhood $\tilde{U}$ of $\Delta$ in $U$ and a regular function
$h_s:\tilde{U}\rightarrow\mathbb{R},$ and a
$\mathcal{C}^{(n+1)k-ns}$ function $g_s:\tilde{U}\rightarrow\mathbb{R}$  such that
$$f-h_s=g_su^s\mbox{ on }\tilde{U}.$$

In order to prove the claim it is sufficient to check by induction on $s$ that Condition(s) holds true for
$s=k$. First using \cite[Lemma 3.6]{bib18} we get an open neighborhood $\tilde{U}$ of $\Delta$ in $U$ and a regular function
$h_1:\tilde{U}\rightarrow\mathbb{R}$ such that
$$h_1|_{\overline{\Gamma}^Z\cap\tilde{U}}=f|_{\overline{\Gamma}^Z\cap\tilde{U}}$$ for every $(d-1)$-dimensional face $\Gamma$
of $\Delta.$ Since every point of $\Delta$ belongs to at most $n$ faces of $\Delta,$ by Fact~\ref{ckdivision}  we obtain a $\C^{(n+1)k-n}$ function $g_1:\tilde{U}\rightarrow\R$ satisfying 
$$g_1=\frac{f-h_1}{u}\mbox{ on }\tilde{U}\setminus u^{-1}(0),$$ and  hence Condition(1) holds.

Now assume that Condition(s) is true for some $1\leq s< k.$ We now prove that for every $(d-1)$-dimensional face $\Gamma$ of
$\Delta$ the function $g_s|_{\overline{\Gamma}^{Z}\cap\tilde{U}}$ is regular. To do this fix $\Gamma.$ We may assume, composing $u, g_s, h_s$ and $f$ with an affine change of coordinates in $\R^n$ if necessary, that   $H_{\Gamma}(x_1,\ldots,x_n)=x_1.$ In particular, $\overline{\Gamma}^Z=\{x_1=0\}$ and
$$\frac{\partial^s}{\partial x_1^s}(f-h_s)|_{\overline{\Gamma}^Z\cap\tilde{U}}=s!g_s(\frac{\partial u}{\partial x_1})^s|_{\overline{\Gamma}^Z\cap\tilde{U}}.$$
The left-hand side of the equation above is (by assumption) a regular function, and $$\{\frac{\partial
u}{\partial x_1}|_{\overline{\Gamma}^Z\cap\tilde{U}}=0\}$$ is the union of hypersurfaces in
$\overline{\Gamma}^Z\cap\tilde{U}$ for which $\frac{\partial u}{\partial x_1}|_{\overline{\Gamma}^Z\cap\tilde{U}}$ is a
minimal defining function. Therefore the left-hand side is divisible by $(\frac{\partial u}{\partial
x_1})^s|_{\overline{\Gamma}^Z\cap\tilde{U}},$ which proves that $g_s|_{\overline{\Gamma}^{Z}\cap\tilde{U}}$ is regular.

Using \cite[Lemma 3.6]{bib18} and shrinking an open neighborhood $\tilde{U}$ of $\Delta$ if necessary, we get a regular
function $w:\tilde{U}\rightarrow\mathbb{R}$ such that
$$w|_{\overline{\Gamma}^Z\cap\tilde{U}}=g_s|_{\overline{\Gamma}^Z\cap\tilde{U}}$$ for every $(d-1)$-dimensional face $\Gamma$
of $\Delta.$ By Fact~\ref{ckdivision}, there is a $\C^{(n+1)k-n(s+1)}$ function
$g_{s+1}:\tilde{U}\rightarrow\R$ satisfying
$${g_{s+1}}=\frac{g_s-w}{u} \mbox{ on }\tilde{U}\setminus u^{-1}(0).$$  In particular,
we have $$g_s=u{g_{s+1}}+w \mbox{ on }\tilde{U}$$ which in conjunction with Condition(s) gives
$$f-h_s=(u{g_{s+1}}+w)u^s=g_{s+1}u^{s+1}+wu^s \mbox{ on }\tilde{U}.$$ Since $h_s+wu^s$ is a regular function on $\tilde{U}$, Condition(s+1)
follows, which completes the proof of Claim~\ref{claim03}.\qed\vspace*{3mm}

Now we easily complete the proof of Lemma \ref{cubeextension} for $n=d$. Namely, by Claim \ref{claim03} we
can write $f=gu^k+h,$ where $g$ is a $\mathcal{C}^k$ function and $h$ is regular. Then it is sufficient to define
$\tilde{f}:=\tilde{g}u^k+h,$ where $\tilde{g}$ is a polynomial approximation of $g$ in the $\mathcal{C}^k$ topology.

Next, assume that $n>d.$ Applying an affine change of coordinates in $\R^n$ if necessary, we may assume that 
$\Delta\subset\R^d\times\{0\}^{n-d}.$ In particular, $\Delta=\Delta'\times\{0\}^{n-d},$ where $\Delta'$ is a $d$-simplex in $\R^d.$ Denote the variables in $\mathbb{R}^n$ by $(x_1,\ldots,x_d,z_1,\ldots,z_{n-d}),$
$x=(x_1,\ldots,x_d), z=(z_1,\ldots,z_{n-d})$ and put $z^{\beta}:=z_1^{\beta_1}\cdots z_{n-d}^{\beta_{n-d}}$ for $\beta=(\beta_1,\ldots,\beta_{n-d})\in\mathbb{N}^{n-d}.$
We look for a function $\tilde{f}$ of the form
$$\tilde{f}(x, z)=\sum_{|\beta|\leq k}z^{\beta}f_{\beta}(x),$$ where the $f_{\beta}$ are regular functions on an open neighborhood of $\Delta'$ in $\R^d$
such that
$\partial^{(\alpha,\beta)}\tilde{f}$ is close to $\partial^{(\alpha,\beta)}f$ on $\Delta=\Delta'\times\{0\}^{n-d}$  and the interpolation condition is satisfied for all $(d-1)$-dimensional faces of $\Delta,$ for suitable $(\alpha,\beta)\in\mathbb{N}^d\times\mathbb{N}^{n-d}.$

Note that for $\tilde{f}$ as above and for each $(\alpha,\beta)\in\mathbb{N}^d\times\mathbb{N}^{n-d}$ we have
$$\partial^{(\alpha,\beta)}\tilde{f}(x,0\ldots,0)=c_{\beta}\partial^{\alpha}f_{\beta}(x),$$ where
for $\beta=(\beta_1,\ldots,\beta_{n-d}),$ it holds $ c_{\beta}=\beta_1!\ldots\beta_{n-d}!\neq 0.$
Hence, for every $\beta\in\mathbb{N}^{n-d}$ with $|\beta|\leq k,$ it is sufficient to find a regular function $f_{\beta}$ on an open neighborhood $\tilde{V}$ of $\Delta'$ in $\R^d$ such that:
$$c_{\beta}f_{\beta}\mbox{ approximates }\partial^{(0,\ldots,0,\beta)}f(\cdot,0,\ldots,0)\mbox{ on }\Delta' \mbox{ in the }\mathcal{C}^{k-|\beta|}\mbox{ topology }$$ and 
$$\partial^{\alpha}c_{\beta}f_{\beta}|_{\overline{\Gamma}^{Z}\cap\tilde{V}}=\partial^{(\alpha,\beta)}f(\cdot,0,\ldots,0)|_{\overline{\Gamma}^{Z}\cap\tilde{V}}$$ for every $(d-1)$-dimensional face $\Gamma$ of $\Delta'$ and every $|\alpha|\leq k-|\beta|-1.$ 

We obtain $f_{\beta}$ required above by applying the already proved version
of Lemma~\ref{cubeextension} for $n=d$ with $\Delta$ replaced by $\Delta'$ and $k$ replaced by $k-|\beta|,$ and $f$ replaced by the 
$\mathcal{C}^{(k-|\beta|)(d+1)}$ function $\partial^{(0,\ldots,0,\beta)}f(\cdot,0,\ldots,0).$ The proof is complete.\qed\vspace*{2mm}

The final goal of this section is a variant of Lemma~\ref{cubeextension} where the simplex $\Delta$ is replaced by an algebraic cell. As a preliminary step, the following lemma will be useful.

\begin{lemma}\label{lemma-4-3}
Let $\varphi:C\rightarrow\Delta$ be a biregular isomorphism between an algebraic cell $C$ in $\R^n$ and a $d$-simplex $\Delta$ in $\R^n.$ Then there is a regular bijection $F:U\rightarrow V$ between open neighborhoods $U, V$ of $C, \Delta,$ respectively, in $\R^n$ such that the following hold:\vspace*{2mm}\\
(t) $F|_{U\cap\overline{C}^Z}:U\cap\overline{C}^Z\rightarrow V\cap\overline{\Delta}^Z$ is a biregular isomorphism,\\
(u) $F|_C=\varphi,$\\
(v) the inverse of $F$ is a map of class $\mathcal{C}^{\infty}$,\\
(w) the restrictions of all partial derivatives (of any order) of $F^{-1}$ to $\overline{\Delta}^Z\cap V$ are regular maps.  
\end{lemma}
\proof First note that, by the definition of regular map, $\varphi$ can be extended to a biregular isomorphism $\tilde{\varphi}:X\rightarrow Y,$ where $X, Y$ are Zariski open neighborhoods of $C, \Delta$ in $\overline{C}^Z, \overline{\Delta}^Z\subset\R^n,$ respectively. Note that $Y$ is a nonsingular variety, hence $X$ is also such.

Let $A, B$ be Euclidean open contractible neighborhoods of $C, \Delta$ in $X, Y,$ respectively, such that $\tilde{\varphi}$ induces a $\C^{\infty}$ diffeomorphism between $A$ and $B.$ Then the normal bundles of $A$ and $B$ in $\R^n$ are trivial. Therefore, using tubular neighborhoods of $A$ and $B$ in $\R^n,$ we can choose
Euclidean open neighborhoods $U, V$ of $C, \Delta$ in $\R^n$ and a $\mathcal{C}^{\infty}$ diffeomorphism $\hat{\varphi}: U\rightarrow V$ such that $U\cap X=U\cap\overline{C}^Z,$
$\hat{\varphi}|_{U\cap X}=\tilde{\varphi}|_{U\cap X}$ and $\hat{\varphi}(U\cap X)=V\cap Y=V\cap\overline{\Delta}^Z.$

Now $\hat{\varphi}$ can be approximated, after shrinking $U, V$ if necessary,  by a regular map $F:U\rightarrow V$ such that $F|_{U\cap X}=\hat{\varphi}|_{U\cap X}$ and $F$ is a diffeomorphism of class $\mathcal{C}^{\infty}$. 
Indeed, first note that $X$ can be expressed as $X=W\cap\overline{C}^Z,$ where $W$ is a Zariski open subset of $\R^n.$ It follows that the regular map $\tilde{\varphi}:X\rightarrow Y\subset\R^n$ can be extended to a regular map $\xi:W\rightarrow\R^n.$ Now, let $\mu_1,\ldots,\mu_s$ be regular functions that generate the ideal of all regular functions on $W$ vanishing on $X.$ Shrinking $U$ if necessary, we may assume that $U$ is a subset of $W.$ Let $\lambda:U\rightarrow\R$ be a $\C^{\infty}$ function vanishing on $U\cap X.$ Since $X$ is nonsingular, for every point $x\in U\cap X,$ the germ of $\lambda$ at $x$ can be written as a linear combination of the germs of the $\mu_i$ at $x$ with coefficients that are $\C^{\infty}$
function-germs $(U,x)\rightarrow\R.$ Using partition of unity, we can express $\lambda$ as
$$\lambda=\lambda_1\mu_1|_U+\cdots+\lambda_s\mu_s|_U,$$ where $\lambda_1,\ldots,\lambda_s$ are $\C^{\infty}$ functions on $U.$
By Weierstrass' theorem, the functions $\lambda_i$ can be approximated in the compact-open $\C^1$ topology by regular functions on $U,$ so $\lambda$ can be approximated by regular functions on $U$ vanishing on $U\cap X.$ Thus, we can find a regular map $\eta: U\rightarrow\R^n$ which vanishes on $U\cap X$ and is arbitrarily close in the compact-open $\C^1$ topology to the $\C^{\infty}$ map $\hat{\varphi}-\xi|_U:U\rightarrow\R^n.$ In view of    
\cite[Lemma 1.3, p. 36]{Hi}, it is sufficient to
define $F:U\rightarrow V$ (after shrinking $U$ and $V$ if necessary) by the formula 
$$F(x):=\xi(x)+\eta(x) \mbox{ for all }x\in U.$$

It remains to show that $F$ satisfies  (w). Let $x=(x_1,\ldots,x_n),$ $y=(y_1,\ldots,y_n)=F(x).$
\begin{claim}\label{claim-4-4}For every $\mathcal{C}^{\infty}$ function $g:V\rightarrow\R$ and every regular function $h:U\rightarrow\R$ such that $$g(F(x))=h(x)\mbox{ for all }x\in U,$$ the following holds: $\frac{\partial g}{\partial y_i}|_{V\cap Y}$ is a regular function for $i=1,\ldots,n.$  
\end{claim}
\noindent\textit{Proof of Claim~\ref{claim-4-4}.} We have $$d_{F(x)}g\circ d_xF=d_xh.$$ The matrix of $d_xF$ is invertible and has entries which are regular functions in $x.$ Hence the matrix of $d_xh\circ(d_xF)^{-1}$ is a row whose entries are regular functions in $x.$  It follows that $\frac{\partial g}{\partial y_i}(F(x))$ is a regular function, for $i=1,\ldots, n.$ Since $F^{-1}|_{V\cap Y}=(\hat{\varphi}|_{U\cap X})^{-1}$ is a regular function, we obtain that $\frac{\partial g}{\partial y_i}(y)=\frac{\partial g}{\partial y_i}(F(F^{-1}(y)))$ is regular on $V\cap Y,$ which proves the claim. \qed\vspace*{2mm} 

Set $G:=F^{-1},$ $G(y)=(g_1(y),\ldots,g_n(y)).$ We will check that for every $j=1,\ldots,n$ and every $\alpha\in\mathbb{N}^n,$ there is a regular function $h_{j,\alpha}$ such that
$$\partial^{\alpha}g_j(F(x))=h_{j,\alpha}(x)\mbox{  for all }x\in U.$$ Once this is done, the proof of Lemma~\ref{lemma-4-3} will be complete by Claim~\ref{claim-4-4}. We proceed by induction on $|\alpha|.$ For $|\alpha|=0$ the assertion is obvious since $g_j(F(x))=x_j,$ for $j=1,\ldots,n.$ Fix any $\alpha\in\mathbb{N}^n$ and $j=1,\ldots,n.$ Put $g:=\partial^{\alpha}g_j.$ By the induction hypothesis, there is a regular function $h_{j,\alpha}$ such that
$$g(F(x))=h_{j,\alpha}(x),\mbox{ for all }x\in U.$$ Therefore, $d_{F(x)}g=d_xh_{j,\alpha}\circ (d_x F)^{-1}.$ As in the proof of Claim~\ref{claim-4-4}, we observe that the matrix of $d_xh_{j,\alpha}\circ (d_x F)^{-1}$ is a row whose entries are regular functions,  hence $\frac{\partial g}{\partial y_i}(F(x))$ is a regular function for every $i=1,\ldots,n,$ which completes the proof of Lemma~\ref{lemma-4-3}.\qed\vspace*{2mm}

%


Let $\mathrm{dist}$ be the metric on $\R^n$ induced by the Euclidean norm $||\cdot||.$ Given  a point $a\in \R^n$ and a constant $\varepsilon>0,$ we denote by $B(a,\varepsilon)$ the Euclidean open ball in $\R^n$ with center at $a$ and radius $\varepsilon.$ For $\alpha=(\alpha_1,\ldots,\alpha_n)\in\mathbb{N}^n$ and $u=(u_1,\ldots,u_n)\in\R^n,$ we set $\alpha!:=\alpha_1!\cdots\alpha_n!$ and $u^{\alpha}:=u_1^{\alpha_1}\cdots u_n^{\alpha_n}.$


Let us recall two well known facts in the form suitable for further applications.
\begin{theorem}[Taylor expansion]\label{Taylor} Let $k$ be a nonnegative integer and let $f:\Omega\rightarrow\mathbb{R}$ be a function of class $\mathcal{C}^{k+1}$ defined on an open set $\Omega$ in $\mathbb{R}^n.$ Let $\Omega_0$ be an open relatively compact subset of $\Omega.$ Then there is a constant
$M\geq 0$ such that
$$|f(x)-\sum_{|\alpha|\leq k}\frac{\partial^{\alpha}f(a)}{\alpha !}(x-a)^{\alpha}|\leq
M||x-a||^{k+1},$$ for all $x,a\in\Omega_0$ with $B(a,||x-a||)\subset\Omega_0.$\end{theorem}

For the first part of the following result see \cite[Th\'eor\`eme 1, p. 85]{Loj}. The special case regarding semilinear sets is  simple and can be obtained by a direct argument.

\begin{theorem}[Regularly situated semialgebraic sets]\label{separation} Let $A, B$ be closed semialgebraic subsets of
$\mathbb{R}^n$ with $A\cap B\neq \emptyset$ and let $\Omega$ be an open bounded subset of $\R^n.$ Then there are constants $\varepsilon>0, r\geq 1$ such that 
$$\mathrm{dist}(x,A)+\mathrm{dist}(x,B)\geq\varepsilon{ }\mathrm{dist}(x,A\cap B)^r \mbox{ for all } x\in\Omega.$$
Moreover, if  $A, B$ are closed semilinear sets, then the inequality above holds with $r=1$ and some constant $\varepsilon>0.$
\end{theorem}

The following proposition will be useful.
\begin{proposition}\label{WhExt} Let $\mathcal{F}$ be a finite family of compact semialgebraic subsets
of $\mathbb{R}^n$ and let $k$ be a nonnegative integer. Then there is a positive integer $p=p(\mathcal{F},k)$ such that the following holds. For every family $\{f_K:U_K\rightarrow\R\}_{K\in\mathcal{F}}$ of $\mathcal{C}^{p+1}$ functions, where each $U_K$ is an open neighborhood of $K$ in $\R^n,$ such that
$$\partial^{\alpha}f_K|_{K\cap L}=\partial^{\alpha}f_L|_{K\cap L}$$ for all $K, L\in\mathcal{F}$ and all $\alpha\in\mathbb{N}^n$ with
$|\alpha|\leq p,$ there is a $\C^k$ function $f: \R^n\rightarrow\R$
satisfying $$\partial^{\alpha} f|_K=\partial^{\alpha} f_K|_K$$ for all $K\in\mathcal{F}$ and all
$\alpha\in\mathbb{N}^n$ with $|\alpha|\leq k.$
\end{proposition}
\proof Let $\Omega$ be a bounded open neighborhood of $\bigcup\mathcal{F}$ in $\R^n.$ For any $K, L\in\mathcal{F}$ with nonempty intersection, let $\varepsilon_{K,L}>0, r_{K,L}\geq 1$ be constants such that 
$$\mathrm{dist}(a,K)+\mathrm{dist}(a,L)\geq \varepsilon_{K,L}\mathrm{dist}(a,K\cap L)^{r_{K,L}}   \mbox{ for all } a\in\Omega.$$
Put $$r:=\max\{r_{K,L}: K, L\in\mathcal{F}, K\cap L\neq\emptyset\}.$$ 
Let $p$ be any integer bounded from below by
$r(k+1).$

Let $\{f_K:U_K\rightarrow\R\}_{K\in\mathcal{F}}$ be a family of $\mathcal{C}^{p+1}$ functions satisfying the assumptions of the proposition. Without loss of generality we assume that each $U_K$ is contained in $\Omega.$ 

Define continuous functions $\tilde{f}_{\alpha}:\bigcup{\mathcal{F}}\rightarrow\R,$ for $\alpha\in\mathbb{N}^n, |\alpha|\leq k,$ as
follows. For every $b\in\bigcup\mathcal{F}$ and $K\in\mathcal{F}$ with $b\in K$ put
$$\tilde{f}_{\alpha}(b)=\partial^{\alpha}f_K(b).$$ Note that $\tilde{f}_{\alpha}(b)$ is well defined (that is,
independently of the choice of $K$ containing $b,$ by assumption).

Now it is sufficient to show that the collection $\{\tilde{f}_{\alpha}\}_{|\alpha|\leq k}$ satisfies the hypotheses of the Whitney
extension theorem \cite[Theorem I, p. 65]{Whit}. More precisely, we prove that for every $c\in\bigcup\mathcal{F}$ we have
$$\hspace*{-21mm}(\ast)\hspace*{21mm}\tilde{f}_{\alpha}(x)=\sum_{|\beta|\leq k-|\alpha|}\frac{\tilde{f}_{\alpha+\beta}(a)}{\beta !}(x-a)^{\beta}+o(||x-a||^{k-|\alpha|}) \mbox{ as }a,x\rightarrow c.$$
It is sufficient to consider the case where $x\in K,$ $a\in L$ and $c\in K\cap L$ for any fixed $K, L\in\mathcal{F},$ and
$x, a$ belong to a ball $B(c,s)$ of radius $s<1/3$
such that  $\overline{B(c,3s)}\subset U_K\cap U_L.$ 

By the Taylor expansion, for $|\alpha|\leq k,$ we have
$$\hspace*{-34.5mm}(+)\hspace*{34.5mm}\partial^{\alpha}{f}_{K}(x)=\sum_{|\beta|\leq k-|\alpha|}\frac{\partial^{\alpha+\beta}{f_K}(a)}{\beta !}(x-a)^{\beta}
+w(x,a),$$ where $|w(x,a)|\leq C(K,\alpha)||x-a||^{k+1-|\alpha|},$ and $C(K,\alpha)$ is a constant obtained by applying Theorem~\ref{Taylor} with $f=\partial^{\alpha}f_K$ and $\Omega_0=B(c,3s).$ In order to prove $(\ast),$ we show that if in the right-hand side of $(+)$ the term ${\partial^{\alpha+\beta}f_K(a)}$ is replaced by $\partial^{\alpha+\beta}f_L(a),$ then $(+)$ remains true with perhaps different $w$. Let $b\in K\cap L$ such that $||a-b||=\mathrm{dist}(a,K\cap L).$ Again by  Theorem~\ref{Taylor} with $\Omega_0=B(c,3s)$, for $|\alpha|+|\beta|\leq k,$ we get
$$\partial^{\alpha+\beta}{f}_{K}(a)=\sum_{|\gamma|\leq p-|\alpha|-|\beta|}\frac{\partial^{\alpha+\beta+\gamma}{f_K}(b)}{\gamma !}(a-b)^{\gamma}
+u(a,b),$$

$$\partial^{\alpha+\beta}{f}_{L}(a)=\sum_{|\gamma|\leq p-|\alpha|-|\beta|}\frac{\partial^{\alpha+\beta+\gamma}{f_L}(b)}{\gamma !}(a-b)^{\gamma}
+v(a,b),$$
where $$|u(a,b)|\leq C(K,\alpha+\beta)||a-b||^{p+1-|\alpha|-|\beta|},$$  $$|v(a,b)|\leq C(L,\alpha+\beta)||a-b||^{p+1-|\alpha|-|\beta|}.$$
In view of the assumption
$$\partial^{\alpha+\beta+\gamma}{f_K}(b)=\partial^{\alpha+\beta+\gamma}{f_L}(b)$$ for $|\alpha|+|\beta|+|\gamma|\leq p,$ we get
$$|\partial^{\alpha+\beta}{f}_{K}(a)-\partial^{\alpha+\beta}{f}_{L}(a)|\leq (C(K,\alpha+\beta)+C(L,\alpha+\beta))||a-b||^{p+1-|\alpha|-|\beta|}.$$ On the other hand, since $K$ and $L$ are closed semialgebraic subsets of $\R^n$ and $a\in L,$ Theorem \ref{separation} implies
$$||a-b||=\mathrm{dist}(a,K\cap L)\leq (\varepsilon_{K,L}^{-1}\mathrm{dist}(a,K))^{\frac{1}{r_{K,L}}}\leq (\varepsilon_{K,L}^{-1}||a-x||)^{\frac{1}{r_{K,L}}}.$$ Therefore 
$$|\partial^{\alpha+\beta}{f}_{K}(a)-\partial^{\alpha+\beta}{f}_{L}(a)|\leq \tilde{C}||a-x||^\frac{p+1-|\alpha|-|\beta|}{r_{K,L}}
\leq \tilde{C}||a-x||^\frac{p+1-|\alpha|-|\beta|}{r}\leq \tilde{C}||a-x||^{k+1-|\alpha|-|\beta|},$$ where $$\tilde{C}=(C(K,\alpha+\beta)+C(L,\alpha+\beta))\varepsilon_{K,L}^{\frac{-p-1+|\alpha|+|\beta|}{r_{K,L}}}.$$

Combining the last conclusion with $(+)$ we have

$$\partial^{\alpha}{f}_{K}(x)=\sum_{|\beta|\leq k-|\alpha|}\frac{\partial^{\alpha+\beta}{f_L}(a)}{\beta !}(x-a)^{\beta}
+o(||x-a||^{k-|\alpha|})\mbox{ as }a,x\rightarrow c.$$
Finally, invoking the definition of $\tilde{f}_{\alpha}$ we get $(\ast),$ which completes the proof.\qed\vspace*{2mm}
\begin{remark}\label{WhRemark} Note that if $\mathcal{F}$ in Proposition~\ref{WhExt} consists of semilinear sets, then as $p=p(\mathcal{F},k)$ we can take
$p=k+1.$ Indeed, by the first paragraph of the proof of Proposition~\ref{WhExt}, $p$ is taken to be any integer satisfying $p\geq r(k+1).$ Moreover,
by Theorem~\ref{separation}, for semilinear sets we have $r=1.$ 
\end{remark}

Given a $d$-dimensional algebraic cell $C$ in $\R^n,$ we denote by $\mathcal{L}_C$ the collection of all $(d-1)$-dimensional faces of $C$.

What follows is the main result of this section.
\begin{lemma}
 \label{mainofthis} Let $C$ be a $d$-dimensional algebraic cell in $\R^n,$ with $d\geq 1.$ Let $l, r$ be nonnegative integers satisfying $l\leq r+1.$ Let $f:U\rightarrow\R$ be a $\C^l$ function defined on an open neighborhood $U\subset\R^n$ of $C.$ Let $\varepsilon>0$ be a constant.
 
 For these data, there exists a constant $\delta>0$ with the following property. If $\{f_K:U_K\rightarrow\R\}_{K\in\mathcal{L}_C}$ is a collection of regular functions, where $U_K\subset U$ is an open neighborhood of $K,$ satisfying\vspace*{2mm}\\
 $(a_{\mathcal{L}_C})$ for every $K\in\mathcal{L}_C$ the functions $f_K$ and $f$ are $\delta$-close on $K$ in the $\C^l$ topology,\\
 $(b_{\mathcal{L}_C})$ $\partial^{\alpha}f_K|_{K\cap L}=\partial^{\alpha}f_L|_{K\cap L}$ for all $K, L\in\mathcal{L}_C$ and all $\alpha\in\mathbb{N}^n$\\ 
 \hspace*{9mm} with $|\alpha|\leq s:=(d+1)(r+1)+2,$\vspace*{2mm}\\
then there exists a regular function $f_C:U_C\rightarrow\R$ defined on an open neighborhood $U_C\subset U$ of $C$ such that\vspace*{2mm}\\
$(a_C)$ $f_C$ and $f$ are $\varepsilon$-close on $C$ in the $\C^l$ topology,\\
$(b_C)$ $\partial^{\alpha}f_C|_K=\partial^{\alpha}f_K|_K$ for all $K\in\mathcal{L}_C$ and all $\alpha\in\mathbb{N}^n$ with $|\alpha|\leq r.$
\end{lemma}
\proof Fix a biregular isomorphism $\varphi: C\rightarrow\Delta$ onto a $d$-simplex
$\Delta$ in $\R^n.$ Shrinking $U$ if necessary, we obtain an open neighborhood $V\subset\R^n$ of $\Delta$ and a regular map $F:U\rightarrow V$ satisfying conditions $(t), (u), (v), (w)$  
in Lemma~\ref{lemma-4-3}.

Choose a constant $\eta>0$ with the following property:\vspace*{2mm}\\
$(C)$ If $\lambda:\tilde{V}\rightarrow\R$ is a $\mathcal{C}^{\infty}$ function defined on an open neighborhood $\tilde{V}\subset V$ of $\Delta$ such 

\hspace*{1.5mm}that the functions $\lambda$ and $f\circ F^{-1}$ are $2\eta$-close on $\Delta$ in the $\C^l$ topology, then the 

\hspace*{1.5mm}functions $\lambda\circ(F|_{F^{-1}(\tilde{V})})$ and $f$ are $\varepsilon$-close on $C$ in the $\C^l$ topology.\vspace*{2mm}

Now choose a constant $\delta>0$ for which the following holds:\vspace*{2mm}\\
$(\mathcal{L}_C)$ For each $K\in\mathcal{L}_{C},$ if $h_K:\Omega_K\rightarrow\R$
is a $\C^{\infty}$ function defined on an open neighborhood 

\hspace*{1.5mm}$\Omega_K\subset U$ of $K$ such that the functions $h_K$ and $f$ are $\delta$-close on $K$ in the $\C^l$ topology, 

\hspace*{1.5mm}then the functions $h_K\circ(F^{-1}|_{F(\Omega_K)})$ and $f\circ F^{-1}$ are $\eta$-close on $F(K)$ in the $\C^l$ topology.\vspace*{2mm}

For $\delta$ as above, let $\{f_K:U_K\rightarrow\R\}_{K\in\mathcal{L}_C}$ be a collection of regular functions satisfying $(a_{\mathcal{L}_C})$ and $(b_{\mathcal{L}_C}).$ Set $V_K:=F(U_K)$ for $K\in\mathcal{L}_C.$ Let $S_K$ be a compact convex semilinear neighborhood of the $(d-1)$-simplex $F(K)$ in $\overline{F(K)}^{Z}$ such that $S_K\subset V_K.$ For each $K\in\mathcal{L}_C$ define a $\C^{\infty}$ function $g_K:V_K\rightarrow \R$ by the formula $$g_K:=f_K\circ(F^{-1}|_{V_K}).$$ According to $(\mathcal{L}_C),$ for every $K\in\mathcal{L}_C$ the functions $g_K$ and $f\circ F^{-1}$ are $\eta$-close on $F(K)$ in the $\C^l$ topology. In view of $(b_{\mathcal{L}_C}),$ we get 
$$\partial^{\alpha}g_K|_{F(K)\cap F(L)}=\partial^{\alpha}g_L|_{F(K)\cap F(L)} $$
for all $K, L\in\mathcal{L}_C$ and all $\alpha\in\mathbb{N}^n$ with $|\alpha|\leq s.$ Moreover, by condition $(w)$ in Lemma~\ref{lemma-4-3}, 
$\partial^{\alpha}g_K|_{V_K\cap\overline{F(K)}^Z}, \partial^{\alpha}g_L|_{V_L\cap \overline{F(L)}^Z}$ are regular functions. Thus, we get
$$\partial^{\alpha}g_K|_{S_K\cap S_L}=\partial^{\alpha}g_L|_{S_K\cap S_L}$$
because $S_K\cap S_L=\emptyset$ if $d=1,$ while for $d\geq 2$ the intersection $S_K\cap S_L$ is a convex semilinear set that contains the $(d-2)$-simplex $F(K)\cap F(L)$ and is contained in the affine-linear $(d-2)$-plane $\overline{F(K)}^Z\cap\overline{F(L)}^Z=\overline{F(K)\cap F(L)}^Z.$

Applying Proposition~\ref{WhExt} and Remark~\ref{WhRemark} to $\mathcal{F}:=\{S_K\}_{K\in\mathcal{L}_C}$ and $\{g_{S_K}:=g_K\}_{S_K\in\mathcal{F}},$ we obtain a $\C^{(d+1)(r+1)}$ function $g:\R^n\rightarrow\R$ satisfying
$$\partial^{\alpha}g|_{S_K}=\partial^{\alpha}g_K|_{S_K}$$
for all $K\in\mathcal{L}_C$ and all $\alpha\in\mathbb{N}^n$ with $|\alpha|\leq (d+1)(r+1).$ By construction, the functions $g$ and $f\circ F^{-1}$ are $\eta$-close on $S$, where $S$ is the union of the family $\mathcal{F},$ in the $\C^l$ topology. Therefore, using $\C^{\infty}$ partition of unity, we can modify $g$ outside a small neighborhood of $S$ in $V,$ so that $g$ and $f\circ F^{-1}$ are in fact $\eta$-close on $\Delta$ in the $\C^l$ topology. Now choose an open neighborhood $W\subset V$ of $\Delta$ such that 
$$S_K\cap W=\overline{F(K)}^Z\cap W\subset \overline{F(K)}^Z\cap V_K$$ for all $K\in\mathcal{L}_C.$ Thus, $\partial^{\alpha}g|_{\overline{F(K)}^Z\cap W}$ is a regular function for all $K\in\mathcal{L}_C$ and all $\alpha\in\mathbb{N}^n$ with $|\alpha|\leq (d+1)(r+1).$ By Lemma~\ref{cubeextension}, there exists a regular function $\tilde{g}:\tilde{V}\rightarrow\R$ defined on an open neighborhood $\tilde{V}\subset W$ of $\Delta$ such that the functions $\tilde{g}$ and $g$ are $\eta$-close on $\Delta$ in the $\C^l$ topology (as $l\leq r+1$), and $$\partial^{\alpha}\tilde{g}|_{\overline{F(K)}^Z\cap\tilde{V}}=\partial^{\alpha}g|_{\overline{F(K)}^Z\cap\tilde{V}}$$ for all $K\in\mathcal{L}_{C}$ and all $\alpha\in\mathbb{N}^n$ with $|\alpha|\leq r.$ Note that the functions $\tilde{g}$ and $f\circ F^{-1}$ are $2\eta$-close on $\Delta$ in the $\C^l$ topology.

Set $U_C:=F^{-1}(\tilde{V})$ and define a regular function by the formula $$f_C:=\tilde{g}\circ(F|_{U_C}).$$ In view of $(C),$ the functions $f_C$ and $f$ are $\varepsilon$-close on $C$ in the $\C^l$ topology. We also have
$$\partial^{\alpha}f_C|_K=\partial^{\alpha}f_K|_{K}$$ for all $K\in\mathcal{L}_C$ and all $\alpha\in\mathbb{N}^n$ with $|\alpha|\leq r.$ Thus, conditions $(a_C)$ and $(b_C)$ hold, which completes the proof.\qed

\section{Properties of uniformly retract rational real\\ algebraic varieties}\label{sec-4}
The proof of Theorem~\ref{mainr} given in this section depends on some properties of uniformly retract rational varieties. These properties, recalled as Definition~\ref{def-4-1} and Theorems~\ref{th-4-2} and \ref{th-4-3} are extracted from Banecki's paper \cite{Ban}

\begin{df}\label{def-4-1}(see \cite[Definition 2.2]{Ban}.) Let $Y$ be a real algebraic variety. A \textit{strong dominating spray} $(q,M,N,\sigma,\tau)$ for $Y$ consists of a nonnegative integer $q$, a Zariski open set $M\subset Y\times\R^q$ containing $Y\times\{0\},$ a Zariski open set $N\subset Y\times Y$ containing the diagonal of $Y\times Y,$ and two regular maps $$\sigma: M\rightarrow Y, \tau:N\rightarrow\R^q$$ subject to the following conditions:\vspace*{2mm}\\
$(i)$ \hspace*{1mm}$(y,\tau(y,z))\in M$ and $\sigma(y,\tau(y,z))=z$ for all $(y,z)\in N,$\\
$(ii)$ $\tau(y,y)=0$ for all $y\in Y.$
\end{df}

Note that condition $(i), (ii)$ yield $$\sigma(y,0)=y\mbox{ for all }y\in Y.$$

\begin{theorem}\label{th-4-2}(see \cite[Observation 1.4, Theorem 2.6]{Ban}.) Every uniformly retract rational real algebraic variety is nonsingular and admits a strong dominating spray. 
\end{theorem}

\begin{theorem}\label{th-4-3}(see \cite[Theorem 1.5]{Ban}.) Let $Y$ be a uniformly retract rational real algebraic variety. Let $X$ be a nonsingular real algebraic variety and let $f: U\rightarrow Y$ be  a $\C^{\infty}$ map defined on an open subset $U$ of $X.$ Assume that $f$ is homotopic to a regular map $U\rightarrow Y.$ Then for every open relatively compact subset $U_0$ of $U$, the restriction $f|_{U_0}$ can be approximated in the compact-open $\C^{\infty}$ topology by regular maps from $U_0$ to~$Y.$
\end{theorem}

The following generalization of Lemma~\ref{mainofthis} is key.
\begin{lemma}\label{lem-4-4}Let $C$ be  a $d$-dimensional algebraic cell in $\R^n,$
with $d\geq 1.$ Let $l, r$ be nonnegative integers satisfying $l\leq r+1.$ Let $Y$ be a uniformly retract rational real algebraic variety which  is a Zariski closed subvariety of $\R^m,$ for some $m.$ Let $f:U\rightarrow Y$ be a $\C^l$ map defined on an open neighborhood $U\subset \R^n$ of $C.$ Let $\varepsilon>0$ be a constant.

For these data, there exists a constant $\delta>0$ with the following property. If $\{f_K:U_K\rightarrow Y\}_{K\in\mathcal{L}_C}$ is a collection of regular maps, where $U_K\subset U$ is an open neighborhood of $K,$ satisfying\vspace*{2mm}\\
$(a_{\mathcal{L}_C,Y})$ for every $K\in\mathcal{L}_C$ the maps $f_K$ and $f$ are $\delta$-close on $K$ in the $\C^l$ topology,\\
$(b_{\mathcal{L}_C,Y})$ $\partial^{\alpha}f_K|_{K\cap L}=\partial^{\alpha}f_L|_{K\cap L}$ for all $K, L\in\mathcal{L}_C$ and all $\alpha\in\mathbb{N}^n$\\ 
\hspace*{14mm}with $|\alpha|\leq s:=(d+1)(r+1)+2,$\vspace*{2mm}\\
then there exists a regular map $f_C:U_C\rightarrow Y$ defined on an open neighborhood $U_C\subset U$ of $C$ such that\vspace*{2mm}\\
$(a_{C,Y})$ $f_C$ and $f$ are $\varepsilon$-close on $C$ in the $\C^l$ topology,\\
$(b_{C,Y})$ $\partial^{\alpha}f_C|_{K}=\partial^{\alpha}f_K|_{K}$ for all $K\in\mathcal{L}_C$ and all $\alpha\in\mathbb{N}^n$ 
with $|\alpha|\leq r.$
\end{lemma}

\proof Since $f$ can be approximated in the compact-open $\C^l$ topology by $\C^{\infty}$ maps, we may assume that $f$ is a $\C^{\infty}$ map.

Shrinking $U$ if necessary, we may assume that $U$ is a contractible set, so $f$ is homotopic to a constant map. Thus by, Theorem~\ref{th-4-3}, $f$ can be approximated in the compact-open $\C^{\infty}$ topology by regular maps from $U$ to $Y.$

According to Theorem~\ref{th-4-2}, there exists a strong dominating spray $(q,M,N,\sigma,\tau)$ for $Y.$ If $g:U\rightarrow Y$ is a regular map sufficiently close to $f,$ then shrinking $U$ again, we get $(g(x),f(x))\in N$ for all $x\in U,$ and hence 
$$h:U\rightarrow\R^q, h(x):=\tau(g(x),f(x))$$ is a well defined $\C^{\infty}$ map.
Moreover, by Definition \ref{def-4-1}(i), $(g(x),h(x))\in M$ and $$f(x)=\sigma(g(x),h(x))\mbox{ for all }x\in U.$$

Let $\{f_K:U_K\rightarrow Y\}_{K\in\mathcal{L}_C}$ be a collection of regular maps satisfying $(a_{\mathcal{L}_C,Y})$ and $(b_{\mathcal{L}_C,Y}).$
If a constant $\delta>0$ is small enough, then for every $K\in\mathcal{L}_C,$ using condition $(a_{\mathcal{L}_C,Y})$ and shrinking $U_K$ if necessary, we get $(g(x),f_K(x))\in N$ for all $x\in U_K,$ so 
$$h_K:U_K\rightarrow\R^q, h_K(x):=\tau(g(x),f_K(x))$$ is a well defined regular map. Additionally, condition $(b_{\mathcal{L}_C,Y})$ implies
$$\partial^{\alpha}h_K|_{K\cap L}=\partial^{\alpha}h_L|_{K\cap L}$$ for all $K, L\in\mathcal{L}_C$ and all $\alpha\in\mathbb{N}^n$ with $|\alpha|\leq s.$ Note that maps $h_K$ and $h$ are as close on $K$ in the $\C^l$ topology as we want, provided that $\delta$ is small enough. Therefore, by Lemma~\ref{mainofthis}, there exists a regular map $h_C:U_C\rightarrow\R^q$ defined on an open neighborhood $U_C\subset U$ of $C$ such that the maps $h_C$ and $h$ are close on $C$ in the $\C^l$ topology, and
$$\partial^{\alpha}h_C|_{K}=\partial^{\alpha}h_K|_{K}$$ for all $K\in\mathcal{L}_C$ and all $\alpha\in\mathbb{N}^n$ with $|\alpha|\leq r.$ Shrinking $U_C$ if necessary, we get $(g(x), h_C(x))\in M$ for all $x\in U_C,$ so 
$$f_C:U_C\rightarrow Y, f_C(x):=\sigma(g(x),h_C(x))$$ is a well defined regular map. If the maps $h_C$ and $h$ are sufficiently close on $C$ in the $\C^l$ topology, then the maps $f_C$ and $f$ are $\varepsilon$-close on $C$ in the $\C^l$-topology. By Definition~\ref{def-4-1}(i), for every $K\in\mathcal{L}_C,$ we get
$$f_K(x)=\sigma(g(x),h_K(x))\mbox{ for all }x\in U_K,$$ and hence
$$\partial^{\alpha}f_C|_K=\partial^{\alpha}f_K|_K$$ for all $\alpha\in\mathbb{N}^n$
with $|\alpha|\leq r.$ Thus, conditions $(a_{C,Y})$ and $(b_{C,Y})$ hold, which completes the proof.\qed\pagebreak

Lemma~\ref{lem-4-4} allows us to prove the following approximation result. 
\begin{proposition}\label{prop-4-5}Let $\mathcal{K}$ be an algebraic complex in $\R^n.$ Let $Y$ be a uniformly retract rational real algebraic variety which is a Zariski closed
subvariety of $\R^m,$ for some $m.$ Let $l$ be a nonnegative integer and let $f:U\rightarrow Y$ be a $\C^l$ map defined on an open neighborhood $U\subset\R^n$ of $|\mathcal{K}|.$ Then for every nonnegative integer $r$ and every constant $\varepsilon>0$ there exists a collection $\{f_K:U_K\rightarrow Y\}_{K\in\mathcal{K}}$ of regular maps, where $U_K\subset U$ is an open neighborhood of $K,$ satisfying\vspace*{2mm}\\
$(a)$ for every $K\in\mathcal{K}$ the maps $f_K$ and $f$ are $\varepsilon$-close on $K$ in the $\C^l$ topology,\\
$(b)$ $\partial^{\alpha}f_K|_{K\cap L}=\partial^{\alpha}f_L|_{K\cap L}$ for all $K,L\in\mathcal{K}$ and all $\alpha\in\mathbb{N}^n$ with $|\alpha|\leq r.$
\end{proposition}
\proof Since $f$ can be approximated in the compact-open $\C^l$ topology by
$\C^{\infty}$ maps, we may assume that $f$ is a $\C^{\infty}$ map.

For a nonnegative integer $d,$ the $d$-skeleton $\mathcal{K}^{(d)}$ of $\mathcal{K}$ is the collection of all algebraic cells in $\mathcal{K}$ of dimension at most $d.$ Clearly, $\mathcal{K}^{(d)}$ is an algebraic complex in $\R^n.$ Using induction on $d,$ we will prove that the proposition holds true with $\mathcal{K}$ replaced by $\mathcal{K}^{(d)},$ for $d\geq 0.$ Once this is done, the proof will be complete since $\mathcal{K}^{(d)}=\mathcal{K}$ for all $d$ sufficiently large. The case of the $0$-skeleton is trivial.

Fix $d\geq 1.$ Let $r\geq l$ be an integer and let $\varepsilon>0$ be a constant. Set $s:=(d+1)(r+1)+2,$ and let $\delta>0$ be a constant that will be determined later. By the induction hypothesis, there exists a collection $\{f_K:U_K\rightarrow Y\}_{K\in\mathcal{K}^{(d-1)}}$ of regular maps, where $U_K\subset U$ is an open neighborhood of $K,$ satisfying\vspace*{2mm}\\
$(\delta)$ for every $K\in\mathcal{K}^{(d-1)}$ the maps $f_K$ and $f$ are $\delta$-close on $K$ in the $\C^l$ topology,\\
$(s)$ $\partial^{\alpha}f_K|_{K\cap L}=\partial^{\alpha}f_L|_{K\cap L}$ for all $K,L\in\mathcal{K}^{(d-1)}$ and all $\alpha\in\mathbb{N}^n$ with $|\alpha|\leq s.$\vspace*{2mm}

According to Lemma~\ref{lem-4-4}, if $\delta>0$ is small enough, then for every $d$-dimensional algebraic cell $C\in\mathcal{K}^{(d)}$ there exists a regular map $f_C:U_C\rightarrow Y$ defined on an open neighborhood $U_C\subset U$ of $C$ such that the maps $f_C$ and $f$ are $\varepsilon$-close on $C$ in the $\C^l$ topology, and $$\partial^{\alpha}f_C|_K=\partial^{\alpha}f_K|_K$$ for all $K\in\mathcal{L}_C$ and all $\alpha\in\mathbb{N}^n$ with $|\alpha|\leq r.$ This completes the proof of the induction step.\qed\vspace*{2mm}

We are now ready to prove the main result of this paper.\vspace*{2mm}\\
\textit{Proof of Theorem~\ref{mainr}.} We may assume that $Y$ is a Zariski closed subvariety of $\R^m,$ for some $m.$ Let $f:U\rightarrow Y$ be a $\C^l$ map defined on an open neighborhood $U\subset\R^n$ of $|\mathcal{K}|.$ Let $p=p(\mathcal{K},k)$ be the positive integer provided by Proposition~\ref{WhExt}. Let $\varepsilon>0$ be a constant. 

Choose a collection $\{f_K:U_K\rightarrow Y\}_{K\in\mathcal{K}}$ of regular maps, where $U_K\subset U$ is an open neighborhood of $K,$ such that the conditions $(a)$ and $(b)$ in Proposition~\ref{prop-4-5} hold with $r$ replaced by $p.$ According to Proposition~\ref{WhExt}, there exists a $\C^k$ map $\tilde{f}:\R^n\rightarrow\R^m$ satisfying $$\partial^{\alpha}\tilde{f}|_K=\partial^{\alpha}{f}_K|_K$$ for all $K\in\mathcal{K}$ and all $\alpha\in\mathbb{N}^n$ with $|\alpha|\leq k.$ By construction, the maps $\tilde{f}$ and $f$ are $\varepsilon$-close
on $|\mathcal{K}|$ in the $\C^l$ topology, and $\tilde{f}(|\mathcal{K}|)\subset Y.$ The proof is complete since $\tilde{f}|_{K}=f_K|_K$ is a regular map for every $K\in\mathcal{K}.$\qed\pagebreak

\noindent\textit{Acknowledgements.} The authors are grateful to Olivier Benoist for useful correspondence. Partial support was provided by the National Science Center (Poland) under grant no. 2022/47/B/ST1/00211.

\end{document}